\documentclass[12pt]{article}%
\usepackage[latin1]{inputenc}
\usepackage{amsfonts,amssymb,latexsym}
\usepackage{amsmath,amsthm}
\usepackage[Sonny]{fncychap}
\usepackage[T1]{fontenc}
\usepackage{pst-all}
\usepackage{pstcol}
\usepackage{indentfirst}
\usepackage[colorlinks=true, bookmarksnumbered=true, bookmarksopen=true,
bookmarksopenlevel=3, pdfstartview=FitH, linkcolor=blue, pdfmenubar=true,
pdftoolbar=true, bookmarks=true,citecolor=green, urlcolor=blue,
filecolor=magenta,plainpages=false,pdfpagelabels,breaklinks]{hyperref}
\usepackage[lmargin=3cm, tmargin=3.5cm,bmargin=3.5cm, rmargin=3cm]{geometry}
\usepackage{amsbsy}
\usepackage{graphicx}
\usepackage{times}
\usepackage{amsmath}
\usepackage{amsfonts}
\usepackage{amssymb}%
\newtheorem{theorem}{Theorem}[section]

\newtheorem{lemma}[theorem]{Lemma}
\newtheorem{remark}[theorem]{Remark}

\numberwithin{equation}{section}
\def\fin {\vskip 0pt \hfill \hbox{\vrule height 5pt width 5pt depth 0pt} \vskip 12pt}
\begin{document}

\title{On the heat equation with nonlinearity and singular anisotropic potential on the boundary}
\author{{Marcelo F. de Almeida}\\{\small Universidade Federal de Sergipe, Departamento de Matem\'{a}tica,} \\{\small {\ CEP 49100-000, Aracaju-SE, Brazil.}}\\{\small \texttt{E-mail:nucaltiado@gmail.com}}\vspace{0.5cm}\\{{Lucas C. F. Ferreira} {\thanks{L. Ferreira was supported by FAPESP and CNPQ,
Brazil. (corresponding author)}}}\\{\small Universidade Estadual de Campinas, Departamento de Matem\'{a}tica,}\\{\small {\ CEP 13083-859, Campinas-SP, Brazil.}}\\{\small \texttt{E-mail:lcff@ime.unicamp.br}}\vspace{0.5cm}\\{Juliana C. Precioso }\\{\small Unesp-IBILCE, Departamento de Matem\'{a}tica,} \\{\small {\ CEP 15054-000, S\~{a}o Jos\'{e} do Rio Preto-SP, Brazil.}}\\{\small \texttt{E-mail:precioso@ibilce.unesp.br}}}
\date{}
\maketitle

\begin{abstract}
This paper concerns with the heat equation in the half-space $\mathbb{R}%
_{+}^{n}$ with nonlinearity and singular potential on the boundary
$\partial\mathbb{R}_{+}^{n}$. We develop a well-posedness theory (\nobreak{without}
using Kato and Hardy inequalities) that allows us to consider critical
potentials with infinite many singularities and anisotropy. Motivated by
potential profiles of interest, the analysis is performed in weak $L^{p}%
$-spaces in which we prove key linear estimates for some boundary operators
arising from the Duhamel integral formulation in $\mathbb{R}_{+}^{n}.$
Moreover, we investigate qualitative properties of solutions like
self-similarity, positivity and symmetry around the axis $\overrightarrow
{Ox_{n}}$. \vspace{0.3cm}

\

\vspace{0.2cm}


\end{abstract}

\setcounter{equation}{0}\setcounter{theorem}{0}

\

\

\section{Introduction}

Heat equations with singular potentials have attracted the interest of many
authors since the work of Baras and Goldstein \cite{Baras1} in the 80's. In a
smooth domain $\Omega\subset\mathbb{R}^{n},$ they studied the Cauchy problem
for the linear heat equation
\begin{equation}
u_{t}-\Delta u-V(x)u=0\text{\ } \label{heat3}%
\end{equation}
with singular potential
\begin{equation}
V(x)=\frac{\lambda}{\left\vert x\right\vert ^{2}} \label{pot-1}%
\end{equation}
and obtained a threshold value $\lambda_{\ast}=\frac{(n-2)^{2}}{4}$ with
$n\geq3$ for existence of positive $L^{2}$-solutions. The potential in
(\ref{pot-1}) is called inverse square (Hardy) potential and is an example of
potential arising from negative power laws. This class of potentials appears
in a number of physical phenomena (see e.g. \cite{FMT1},\cite{FMT2}%
,\cite{Frank},\cite{Landau},\cite{Levy},\cite{Peral-Vazquez} and references
therein) and can be classified according to the number of singularities
(poles), $\sigma$-degree of the singularity (order of the poles), dependence
on directions (anisotropy) and decay at infinity. One of the most difficult
cases is the one of anisotropic critical potentials, namely
\begin{equation}
V(x)=\sum_{i=1}^{l}\frac{v_{i}\left(  \frac{x-x^{i}}{|x-x^{i}|}\right)
}{|x-x^{i}|^{\sigma}}, \label{pot-int}%
\end{equation}
where $v_{i}(z)\in BC(\mathbb{S}^{n-1})$, $x^{i}\in\overline{\Omega},$
$l\in\mathbb{N\cup\{\infty\}},$ and the parameter $\sigma$ is the order of the
poles $\{x^{i}\}_{i=1}^{l}$. The potential is called isotropic (resp.
anisotropic) when the $v_{i}$'s are independent (resp. dependent) of the
directions $\frac{x-x^{i}}{\left\vert x-x^{i}\right\vert }$, that is, they are
constant. In the case $l=1$ (resp. $l>1$)$,$ $V$ is said to be monopolar
(resp. multipolar). The criticality means that $\sigma$ is equal to order of
the PDE inside the domain or of the boundary condition, according to the type
of problem considered. The critical case introduces further difficulties in
the mathematical analysis of the problem because $Vu$ cannot be handled as a
lower order term (see \cite{FMT1}). Examples of (\ref{pot-int}) are
\begin{equation}
V(x)=\sum_{i=1}^{l}\frac{\lambda_{i}}{|x-x^{i}|^{\sigma}}\text{ and }%
V(x)=\sum_{i=1}^{l}\frac{(x-x^{i}).d^{i}}{|x-x^{i}|^{\sigma+1}}\text{,}
\label{pot2}%
\end{equation}
where $x^{i}=(x_{1}^{i},x_{2}^{i},...,x_{n}^{i})$ and $d^{i}\in\mathbb{R}^{n}$
are constant vectors.\textbf{ }In the theory of Schrodinger operators, the
potentials in (\ref{pot2}) are called multipolar Hardy potentials and multiple
dipole-type potentials, respectively (see \cite{FMT1} and \cite{FMT2}).

In this paper, we consider a nonlinear counterpart for (\ref{heat3}) in the
half-space with critical singular boundary potential, which reads as \
\begin{align}
\partial_{t}u  &  =\Delta u\text{ \ in }\Omega,\;t>0\label{h1}\\
\partial_{n}u  &  =h(u)+V(x^{\prime})u\text{ \ \ in }\partial\Omega
,\;t>0\label{h2}\\
u(x,0)  &  =u_{0}(x),\text{ \ in }\Omega, \label{h3}%
\end{align}
where $\Omega=\mathbb{R}_{+}^{n},$ $n\geq3$ and $\partial_{n}=-\partial
_{x_{n}}$ stands for the normal derivative on $\partial\mathbb{R}_{+}^{n}$.
For the nonlinear term, we assume that the function $h:\mathbb{R\rightarrow
R}$ satisfies $h(0)=0$ and%
\begin{equation}
\left\vert h(a)-h(b)\right\vert \leq\eta\left\vert a-b\right\vert \left(
\left\vert a\right\vert ^{\rho-1}+\left\vert b\right\vert ^{\rho-1}\right)  ,
\label{nonli-1}%
\end{equation}
where $\rho>1$ and the constant $\eta$ is independent of $a,b\in\mathbb{R}$. A
classical example of $h$ satisfying these conditions is $h(u)=\pm\left\vert
u\right\vert ^{\rho-1}u.$

Our goal is to develop a global-in-time well-posedness theory for (\ref{h1}%
)-(\ref{h3}), under smallness conditions on certain weak norms of $u_{0},V$,
that allows to consider critical potentials on the boundary with infinite many
singularities. For that matter, we employ the framework of weak-$L^{p}$ spaces
(i.e., $L^{(p,\infty)}$-spaces) and take $V\in L^{(n-1,\infty)}(\partial
\mathbb{R}_{+}^{n}).$ Since $L^{p}(\partial\mathbb{R}_{+}^{n})$ contains only
trivial homogeneous functions, a motivation naturally appears for considering
weak-$L^{p}$ spaces. In fact, due to Chebyshev's inequality, we have the
continuous inclusion $L^{p}(\partial\mathbb{R}_{+}^{n})\subset L^{(p,\infty
)}(\partial\mathbb{R}_{+}^{n})$ and then $L^{(p,\infty)}$ can be regarded as a
natural extension of $L^{p}$ which contains homogeneous functions of degree
$\sigma=$ $-(n-1)/p.$ The critical case for (\ref{h1})-(\ref{h3}) with
potential (\ref{pot-int}) corresponds to $\sigma=1$ (so, $p=n-1$) and we have that%

\begin{equation}
\left\Vert V\right\Vert _{L^{(n-1,\infty)}(\partial\mathbb{R}_{+}^{n})}%
\leq\sum_{i=1}^{l}\sup_{x^{\prime}\in\mathbb{S}^{n-2}}\left\vert
v_{i}(x^{\prime})\right\vert \left\Vert \frac{1}{|x^{\prime}-x^{i}%
|}\right\Vert _{L^{(n-1,\infty)}(\partial\mathbb{R}_{+}^{n})}\leq C\sum
_{i=1}^{l}\sup_{x^{\prime}\in\mathbb{S}^{n-2}}\left\vert v_{i}(x^{\prime
})\right\vert , \label{pot-aux-1}%
\end{equation}
where%
\begin{equation}
C=\left\Vert \,|x^{\prime}|^{-1}\right\Vert _{L^{(n-1,\infty)}(\partial
\mathbb{R}_{+}^{n})}<\infty. \label{aux-weak-1}%
\end{equation}
Of special interest is when the set $\{x^{i}\}_{i=1}^{l}\subset\partial
\mathbb{R}_{+}^{n}$ and so $V$ has a number $l$ of singularities on the
boundary which can be infinite provided that the infinite sum in
(\ref{pot-aux-1}) is finite.

We address (\ref{h1})-(\ref{h3}) by means of the following equivalent integral
formulation
\begin{equation}
u(x,t)=\int_{\mathbb{R}_{+}^{n}}G(x,y,t)u_{0}(y)dy+\int_{0}^{t}\int
_{\partial\mathbb{R}_{+}^{n}}G(x,y^{\prime},t-s)\left[  h(u)+Vu\right]
(y^{\prime},s)dy^{\prime}ds \label{int}%
\end{equation}
where $G(x,y,t)$ is the heat fundamental solution in $\mathbb{R}_{+}^{n}$
given by
\begin{equation}
G(x,y,t)=(4\pi t)^{-\frac{n}{2}}\left[  e^{-\frac{|x-y|^{2}}{4t}}%
+e^{-\frac{|x-y^{\ast}|^{2}}{4t}}\right]  ,\;x,y\in\overline{\mathbb{R}%
_{+}^{n}},\;t>0, \label{Gxyt}%
\end{equation}
with $y^{\ast}=(y^{\prime},-y_{n})$ and $y^{\prime}=(y_{1},\cdots,y_{n-1}%
)\in\partial\mathbb{R}_{+}^{n}$. Here solutions for (\ref{int}) are looked for
in $BC((0,\infty);\mathcal{X}_{p,q})$ where $\mathcal{X}_{p,q}$ is a suitable
Banach space that can be identified with $L^{(p,\infty)}(\mathbb{R}_{+}%
^{n})\times L^{(q,\infty)}(\partial\mathbb{R}_{+}^{n})$. The norm in
$\mathcal{X}_{p,q}$ provides a $L^{(q,\infty)}$-information for $u|_{\partial
\mathbb{R}_{+}^{n}}$ without assuming any positive regularity condition on
$u$. Notice that this space is specially useful in order to treat singular
boundary terms like (\ref{pot-int}). $L^{r}$-versions of $\mathcal{X}_{p,q}$
(i.e. $L^{r_{1}}(\Omega)\times L^{r_{2}}(\partial\Omega)$) was employed in
\cite{Quittner1} and \cite{Fer-Vil1} for studying weak solutions for an
elliptic and parabolic PDE in bounded domains $\Omega,$ respectively. Let us
observe that $\left\Vert |x^{\prime}|^{-1}\right\Vert _{L^{r_{2}}%
(\partial\mathbb{R}_{+}^{n})}=\infty$ for all $1\leq r_{2}\leq\infty$ (compare
with (\ref{aux-weak-1})) which prevents the use of the spaces of
\cite{Fer-Vil1, Quittner1} for our purposes.

Furthermore, we investigate qualitative properties of solutions like
positivity, symmetries (e.g. invariance around the axis $\overrightarrow
{Ox_{n}}$) and self-similarity, under certain conditions on $u_{0},V,h(\cdot
)$. For the latter, the indexes of spaces are chosen so that their norms are
invariant by scaling of (\ref{h1})-(\ref{h2}) (see (\ref{sc1}) below), namely
$p=n(\rho-1)$ and $q=(n-1)(\rho-1).$

Common tools used to handle (\ref{heat3}) and (\ref{h1})-(\ref{h2}) with
critical potentials are the Hardy and Kato inequalities, which read
respectively as
\begin{align}
\frac{(n-2)^{2}}{4}\int_{\mathbb{R}^{n}}\frac{u^{2}}{\left\vert x\right\vert
^{2}}dx  &  \leq\left\Vert \nabla u\right\Vert _{L^{2}(\mathbb{R}^{n})}%
^{2},\text{ }\forall\varphi\in C_{0}^{\infty}(\mathbb{R}^{n}%
)\label{Hardy-Ineq-1}\\
2\frac{\Gamma(\frac{n}{4})^{2}}{\Gamma(\frac{n-2}{4})^{2}}\int_{\partial
\mathbb{R}_{+}^{n}}\frac{u^{2}}{\left\vert x\right\vert }dx  &  \leq\left\Vert
\nabla u\right\Vert _{L^{2}(\mathbb{R}_{+}^{n})}^{2},\text{ }\forall\varphi\in
C_{0}^{\infty}(\overline{\mathbb{R}_{+}^{n}}),\text{ } \label{Hardy-Ineq-2}%
\end{align}
where $\Gamma$ stands for the gamma function. Our approach relies on a
contraction argument in the space $BC((0,\infty);\mathcal{X}_{p,q})$ which
does not require (\ref{Hardy-Ineq-1}) nor (\ref{Hardy-Ineq-2}). For this
purpose, we need to prove estimates in weak-$L^{p}$ for some boundary
operators linked to (\ref{int}). In view of (\ref{pot-int}), these estimates
need to be time-independent and thereby one cannot use time-weighted norms
\textit{ala kato }(see \cite{Kato1} for this type of norm), making things more
difficult-to-treating. This situation leads us to derive boundary estimates in
spirit of the paper \cite{Yamazaki} that dealt with the heat and Stokes
operators inside a half-space (among other smooth domains $\Omega$). So, in a
certain sense, Lemma \ref{prop-key1} can be seen as extensions of Yamazaki's
estimates to boundary operators. Also, before obtaining Lemma \ref{prop-key1},
we need to prove Lemmas \ref{prop-trace} and \ref{nonl-key-lemma} that seems
to have an interest of its own. It is worthy to comment that weak-$L^{p}$
spaces are examples of shift-invariant Banach spaces of local measure for
which global-in-time well-posedness theory of small solutions has been successfully
developed for Navier-Stokes equations (see \cite{Lemarie} for a nice review)
and, more generally, parabolic problems with nonlinearities (and possibly
other terms) defined inside the domain (see \cite{Kozono-Yamazaki}).

Let us review some works concerning heat equations with singular potentials
and nonlinear boundary conditions. The paper of Baras-Goldstein \cite{Baras1}
have motivated many works concerning heat equations with singular potentials.
In these results, Hardy type-inequalities play an important role in both
linear and nonlinear cases. For (\ref{heat3}) (potential defined inside the
domain), we refer the reader to \cite{Cabre-Martel1},\cite{Goldstein2}%
,\cite{Vazquez1} (see also their references) for results on existence,
non-existence, decay and self-similar asymptotic behavior of solutions.
Versions of (\ref{heat3}) with nonlinearities $\pm u^{p}$ and $\pm\left\vert
\nabla u\right\vert ^{p}$ have been studied in \cite{Abdellaoui-Peral2}%
,\cite{Abdellaoui-Peral3}, \cite{Abdellaoui-Peral4},\cite{Chaves}%
,\cite{Karachalios},\cite{Liskevich1},\cite{Reyes1} where the reader can found
results on existence, non-existence, Fujita exponent, self-similarity,
bifurcations, and blow-up. Linear and nonlinear elliptic versions of
(\ref{heat3}) are also often considered in the literature (see e.g.
\cite{Chaves},\cite{FMT1},\cite{FMT2},\cite{Felli1},\cite{Fer-Mesq}%
,\cite{Fer-Mont},\cite{Smets}); as well as the parabolic case, the key tool
used in the analysis is Hardy type inequalities, except by \cite{Fer-Mesq} and
\cite{Fer-Mont}. In these last two references, the authors employed a
contraction argument in a sum of weighted spaces and in a space based on
Fourier transform, respectively. In a bounded domain $\Omega$ and half-space
$\mathbb{R}_{+}^{n},$ the nonlinear problem (\ref{h1})-(\ref{h3}) with
$V\equiv0$ has been studied by several authors over the past two decades; see,
e.g., \cite{Arrieta-1},\cite{Arrieta-Nolosca-1},\cite{Garcia-Peral}%
,\cite{Ishige},\cite{Quittner2},\cite{Rodriguez-Bernal} and their references.
In these works, the reader can find many types of existence and asymptotic
behavior results in the framework of $L^{p}$-spaces. For $V\in L^{\infty
}(\partial\Omega)$ and $\Omega$ a bounded smooth domain, results on
well-posedness and attractors can be found in \cite{Arrieta-Nolasco-2}. The
authors of \cite{Daners-1} considered (\ref{h1})-(\ref{h3}) with $h(u)\equiv0$
(linear case) and showed $L^{p}$-estimates of solutions, still for $V\in
L^{\infty}(\partial\Omega)$ (see also \cite{Daners-2} for the elliptic case).
In \cite{Ishige2}, the authors studied the linear case of (\ref{h1}%
)-(\ref{h3}) in a half-space and considered the singular critical potential
$V(x^{\prime})=\frac{\lambda}{|x^{\prime}|}.$ For compactly supported data
$u_{0}\in C_{0}(\mathbb{R}_{+}^{n}),$ they obtained a threshold value for
existence of positive solutions by using the Kato inequality
(\ref{Hardy-Ineq-2}).

In this paragraph, we summarize the novelties of the present paper in
comparison with the previous ones. Our results provide a global-in-time well-posedness theory
for (\ref{h1})-(\ref{h3}) in a framework that is larger than $L^{p}$-spaces
and seems to be new in the study of parabolic problems with nonlinear boundary
conditions. Also, among others, it allows us to consider critical potentials
on the boundary with infinite many singularities which are not covered by
previous results. As pointed out above, a remarkable difference is that the
approach employed here does not use Hardy nor Kato inequalities, being based
on boundary estimates on weak-$L^{p}$ spaces. Since the smallness condition on
$u_{0}$ is with respect to the weak norm of such spaces, some initial data
with large $L^{p}$ and $H^{s}$-norms can be considered. Results on
self-similarity and axial-symmetry are naturally obtained due to the choice of
the space indexes and to the symmetry features of the linear operators arising
in the integral formulation (\ref{int}).

The plan of this paper is the following. In the next section we summarize some
basic definitions and properties on Lorentz spaces. In section \ref{fs-r} we
define suitable time-functional spaces and state our results, which are proved
in section \ref{proofs}.

\section{Preliminaries}

\label{section2}

In this section we fix some notations and summarize basic properties about
Lorentz spaces that will be used throughout the paper. For further details, we
refer the reader to \cite{BS},\cite{BL}.

For a point $x\in\overline{\mathbb{R}_{+}^{n}},$ we write $x=(x^{\prime}%
,x_{n})$ where $x^{\prime}=(x_{1},x_{2},\ldots,x_{n-1})\in\mathbb{R}^{n-1}$
and $x_{n}\geq0$. The Lebesgue measure in a measurable $\Omega\subset$
$\mathbb{R}^{n}$ will be denoted by either $|\cdot|$ or $dx$. In the case
$\Omega=\mathbb{R}_{+}^{n}$, one can express $dx=dx^{\prime}dx_{n}$ where
$dx^{\prime}$ stands for Lebesgue measure on $\partial\mathbb{R}_{+}%
^{n}=\mathbb{R}^{n-1}$. Given a subset $\Omega\subset\mathbb{R}^{n},$ the
distribution function and rearrangement of a measurable function
$f:\Omega\rightarrow\mathbb{R}$ is defined respectively by
\[
\lambda_{f}(s)=\left\vert \{x\in\Omega:|f(x)|>s\}\right\vert \text{ \ and
}f^{\ast}(t)=\inf\{s>0:\lambda_{f}(s)\leq t\}\text{, }t>0.
\]

The \textit{Lorentz space} $L^{(p,r)}=L^{(p,r)}(\Omega)=L^{(p,r)}%
(\Omega,|\cdot|)$ consists of all measurable functions $f$ in $\Omega$ for
which
\begin{equation}
\left\Vert f\right\Vert _{L^{(p,r)}(\Omega)}^{\ast}=%
\begin{cases}
\left[  \int_{0}^{\infty}\left(  t^{\frac{1}{p}}[f^{\ast}(t)]\right)
^{r}\frac{dt}{t}\right]  ^{\frac{1}{r}}<\infty, & 0<p<\infty,1\leq r<\infty\\
& \\
\displaystyle\sup_{t>0}t^{\frac{1}{p}}[f^{\ast}(t)]<\infty, & 0<p<\infty
,\,r=\infty.
\end{cases}
\label{almost_norm}%
\end{equation}
We have that $L^{p}(\Omega)=L^{(p,p)}(\Omega)$ and $L^{(p,\infty)}$ is also
called weak-$L^{p}$ or Marcinkiewicz space. The quantity (\ref{almost_norm})
is not a norm in $L^{(p,r)}$, however it is a complete quasi-norm. Considering%
\[
f^{\ast\ast}(t)=\frac{1}{t}\int_{0}^{t}f^{\ast}(s)ds,
\]
we can endow $L^{(p,r)}$ with the quantity $\Vert\cdot\Vert_{L^{(p,r)}}$
obtained from (\ref{almost_norm}) with $f^{\ast\ast}$ in place of $f^{\ast}$.
For $1<p\leq\infty,$ we have that $\Vert\cdot\Vert_{(p,r)}^{\ast}\leq
\Vert\cdot\Vert_{(p,r)}\leq\frac{p}{p-1}\Vert\cdot\Vert_{(p,r)}^{\ast}$ which
implies that $\Vert\cdot\Vert_{(p,r)}^{\ast}$ and $\Vert\cdot\Vert_{(p,r)}$
induce the same topology on $L^{(p,r)}$. Moreover, the pair $(L^{(p,r)}%
,\Vert\cdot\Vert_{L^{(p,r)}})$ is a Banach space. From now on, for
$1<p\leq\infty$ we consider $L^{(p,r)}$ endowed with $\Vert\cdot
\Vert_{L^{(p,r)}}$, except when explicitly mentioned.

If, for $\lambda>0,$ $\lambda\Omega=\{\lambda x:x\in\Omega\}$ is the
\textit{dilation} of the domain $\Omega$, then
\begin{equation}
\Vert f(\lambda x)\Vert_{L^{(p,r)}(\Omega)}=\lambda^{-\frac{n}{p}}\Vert
f(x)\Vert_{L^{(p,r)}(\Omega)},\label{scaling-aux-1}%
\end{equation}
provided that $\Omega$ is invariant by dilations, i.e., $\Omega=\lambda\Omega$.

For $1\leq q_{1}\leq p\leq q_{2}\leq\infty$ with $1<p\leq\infty$, the
continuous inclusions hold true%
\[
L^{(p,1)}\subset L^{(p,q_{1})}\subset L^{p}\subset L^{(p,q_{2})}\subset
L^{(p,\infty)}.
\]
The dual space of $L^{(p,r)}$ is $L^{(p^{\prime},r^{\prime})}$ for $1\leq
p,r<\infty$. In particular, the dual of $L^{(p,1)}$ is $L^{(p^{\prime}%
,\infty)}$ for $1\leq p<\infty.$

H\"{o}lder's inequality works well in Lorentz spaces (see \cite{ONeil}).
Precisely, if $1<p_{1},p_{2},p_{3}<\infty$ and $1\leq r_{1},r_{2},r_{3}%
\leq\infty$ with $1/p_{3}=1/p_{1}+1/p_{2}$ and $1/r_{3}\leq1/r_{1}+1/r_{2},$
then
\begin{equation}
\Vert fg\Vert_{L^{(p_{3},r_{3})}}\leq C\Vert f\Vert_{L^{(p_{1},r_{1})}}\Vert
g\Vert_{L^{(p_{2},r_{2})}}, \label{holder}%
\end{equation}
where $C>0$ is a constant independent of $f,g$.

Finally we recall some interpolation property of Lorentz spaces. For
$0<p_{1}<p_{2}<\infty,\ 0<\theta<1,\ \frac{1}{p}=\frac{1-\theta}{p_{1}}%
+\frac{\theta}{p_{2}}$ and $1\leq r_{1},r_{2},r\leq\infty$, we have that (see
\cite[Theorems 5.3.1, 5.3.2]{BL})
\begin{equation}
\left(  L^{(p_{1},r_{1})},L^{(p_{2},r_{2})}\right)  _{\theta,r}=L^{(p,r)},
\label{interp-1}%
\end{equation}
where $(X,Y)_{\theta,r}$ stands for the real interpolation space between $X$
and $Y$ constructed via the $K_{\theta,q}$-method. It is well known that
$(\cdot,\cdot)_{\theta,r}$ is an exact interpolation functor of exponent
$\theta$ on the categories of quasi-normed and normed spaces. When
$0<p_{1}\leq1,$ the property (\ref{interp-1}) should be considered with
$L^{(p_{1},r_{1})}$ endowed with the complete quasi-norm $\left\Vert
\cdot\right\Vert _{L^{(p_{1},r_{1})}}^{\ast}$ instead of $\left\Vert
\cdot\right\Vert _{L^{(p_{1},r_{1})}}.$

\section{Functional setting and results}

\label{fs-r}

Before starting our results, we define suitable function spaces where
(\ref{int}) will be handled. If the potential $V$ is a homogeneous function of
degree $-1$, that is, $V(y)=\lambda V(\lambda y)$ for all $y\in\partial
\mathbb{R}_{+}^{n}$, then $u_{\lambda}(x,t)=\lambda^{\frac{1}{\rho-1}%
}u(\lambda x,\lambda^{2}t)$ verifies (\ref{h1})-(\ref{h2}), for each fixed
$\lambda>0$, provided that $u(x,t)$ is also a solution. It follows that
(\ref{h1})-(\ref{h2}) has the following scaling
\begin{equation}
u(x,t)\rightarrow u_{\lambda}(x,t)=\lambda^{\frac{1}{\rho-1}}u(\lambda
x,\lambda^{2}t),\;\lambda>0.\text{ } \label{sc1}%
\end{equation}
Making $t\rightarrow0^{+}$ in (\ref{sc1}), one obtains
\begin{equation}
u_{0}(x)\rightarrow u_{0,\lambda}(x,0)=\lambda^{\frac{1}{\rho-1}}u_{0}(\lambda
x), \label{sc2}%
\end{equation}
which gives a scaling for the initial data.

Since the potential $V$ and initial data $u_{0}$ are singular, we need to
treat (\ref{int}) in a \linebreak suitable space of functions without any
positive regularity conditions and time decaying. For that matter, let
$\mathcal{A}$ be the set of measurable functions $f:\overline{\mathbb{R}%
_{+}^{n}}\rightarrow\mathbb{R}$ such that $f|_{\mathbb{R}_{+}^{n}}$ and
$f|_{\partial\mathbb{R}_{+}^{n}}$ are measurable with respect to Lebesgue
$\sigma\text{-algebra}$ on $\mathbb{R}_{+}^{n}\text{ and }\mathbb{R}%
^{n-1}=\partial\mathbb{R}_{+}^{n}$, respectively. Consider the equivalence
relation in $\mathcal{A}$: $f\sim g$ if and only if $f=g$ a.e. in
$\mathbb{R}_{+}^{n}$ and $f|_{\partial\mathbb{R}_{+}^{n}}=g|_{\partial
\mathbb{R}_{+}^{n}}$ a.e. in $\partial\mathbb{R}_{+}^{n}$. Given $1\leq
p,q<\infty$, we set $\mathcal{X}_{p,q}$ as the space of all $f\in
\mathcal{A}/\sim$ such that
\begin{equation}
\Vert f\Vert_{\mathcal{X}_{p,q}}=\Vert f\Vert_{L^{(p,\infty)}(\mathbb{R}%
_{+}^{n})}+\Vert f|_{\partial\mathbb{R}_{+}^{n}}\Vert_{L^{(q,\infty)}%
(\partial\mathbb{R}_{+}^{n})}<\infty.\nonumber
\end{equation}
The pair $(\mathcal{X}_{p,q},\Vert\cdot\Vert_{\mathcal{X}_{p,q}})$ is a Banach
space and can be isometrically identified with $L^{(p,\infty)}(\mathbb{R}%
_{+}^{n})\times L^{(q,\infty)}(\partial\mathbb{R}_{+}^{n}).$ For $p=n(\rho-1)$
and $q=(n-1)(\rho-1),$ we have from (\ref{scaling-aux-1}) that%

\begin{equation}
\Vert\lambda^{\frac{1}{\rho-1}}f(\lambda x)\Vert_{\mathcal{X}_{p,q}}%
=\lambda^{\frac{1}{\rho-1}}\lambda^{-\frac{n}{n(\rho-1)}}\Vert f\Vert
_{L^{(p,\infty)}(\mathbb{R}_{+}^{n})}+\lambda^{\frac{1}{\rho-1}}%
\lambda^{-\frac{n-1}{(n-1)(\rho-1)}}\Vert f|_{\partial\mathbb{R}_{+}^{n}}%
\Vert_{L^{(q,\infty)}(\partial\mathbb{R}_{+}^{n})}=\Vert f\Vert_{\mathcal{X}%
_{p,q}},\nonumber
\end{equation}
and then $\mathcal{X}_{p,q}$ is invariant by scaling (\ref{sc2}).

We shall look for solutions in the Banach space $E=BC((0,\infty);\mathcal{X}%
_{p,q})$ endowed with the norm
\begin{equation}
\left\Vert u\right\Vert _{E}=\sup_{t>0}\left\Vert u(\cdot,t)\right\Vert
_{\mathcal{X}_{p,q}}, \label{norm1}%
\end{equation}
which is invariant by scaling (\ref{sc1}).

\subsection{Existence and self-similarity}

In what follows, we state our well-posedness result.

\begin{theorem}
\label{existence} Let $n\geq3,$ $\rho>1$ with $\frac{\rho}{\rho-1}<n-1,$
$p=n(\rho-1)$ and $q=(n-1)(\rho-1).$ Let $h:\mathbb{R\rightarrow R}$ verify
(\ref{nonli-1}) and $h(0)=0$. Suppose that $V\in L^{(n-1,\infty)}%
(\partial\mathbb{R}_{+}^{n})$ and $u_{0}\in L^{(p,\infty)}(\mathbb{R}_{+}%
^{n})$.

\begin{itemize}
\item[(A)] (Existence and uniqueness) There exist $\varepsilon,\delta
_{1},\delta_{2}>0$ such that if $\Vert V\Vert_{L^{(n-1,\infty)}(\partial
\mathbb{R}_{+}^{n})}<\frac{1}{\delta_{1}}$ and $\Vert u_{0}\Vert
_{L^{(p,\infty)}(\mathbb{R}_{+}^{n})}\leq\frac{\varepsilon}{\delta_{2}}$ then
the integral equation (\ref{int}) has a unique solution $u\in BC((0,\infty
);\mathcal{X}_{p,q})$ satisfying $\sup_{t>0}\Vert u(\cdot,t)\Vert
_{\mathcal{X}_{p,q}}\leq\frac{2\varepsilon}{1-\gamma}$ where $\gamma
=\delta_{1}\Vert V\Vert_{L^{(n-1,\infty)}(\partial\mathbb{R}_{+}^{n})}$.
Moreover, $u(\cdot,t)\rightharpoonup u_{0}$ in $\mathcal{S}^{\prime
}(\mathbb{R}_{+}^{n})$ as $t\rightarrow0^{+}.$

\item[(B)] (Continuous dependence) The solution obtained in item (A) depends
continuously on the initial data $u_{0}$ and potential $V$.
\end{itemize}
\end{theorem}

\bigskip

\begin{remark}
\

\begin{itemize}
\item[(A)] The integral solution $u(\cdot,t)\in\mathcal{X}_{p,q},$ for each
$t>0,$ even requiring only $u_{0}\in L^{(p,\infty)}(\mathbb{R}_{+}^{n}).$ It
is a kind of \textquotedblleft parabolic regularizing effect\textquotedblright
in the sense that solutions verifies a property for $t>0\,\ $that is not
necessarily verified by initial data. Here, it comes essentially from the fact
that $u_{1}(x,t)=\int_{\mathbb{R}_{+}^{n}}G(x,y,t)u_{0}(y)dy$ has a trace
well-defined on $\partial\mathbb{R}_{+}^{n}$ and $u_{1}|_{\partial
\mathbb{R}_{+}^{n}}\in L^{(q,\infty)}(\partial\mathbb{R}_{+}^{n}),$ for each
$t>0$, even if the data $u_{0}$ does not have a trace. Moreover, the estimate
(\ref{est-1}) provides a control on the trace just using the norm of $u_{0}$
in $L^{(p,\infty)}(\mathbb{R}_{+}^{n})$.

\item[(B)] The time-continuity of the solution $u(\cdot,t)\in\mathcal{X}%
_{p,q}$ at $t>0$ comes naturally from the uniform continuity of the kernel
$G(x,y,t)$ (\ref{Gxyt}) on $\overline{\mathbb{R}_{+}^{n}}\times\overline
{\mathbb{R}_{+}^{n}}\times\lbrack\delta,\infty),$ for each fixed $\delta>0$.

\item[(C)] A standard argument shows that solutions of (\ref{int}) obtained in
Theorem \ref{existence} verifies (\ref{h1})-(\ref{h3}) in the sense of distributions.
\end{itemize}
\end{remark}

\bigskip

Since the spaces in which we look for solutions are invariant by scaling
(\ref{sc1}), it is natural to ask about existence of self-similar solutions.
This issue is considered in the next theorem.

\begin{theorem}
\label{self-similar}Let $u$ be the mild solution obtained in Theorem
\ref{existence} corresponding to the triple $(u_{0},V,h(\cdot))$. If $u_{0},V$
and $h(\cdot)$ are homogeneous functions of degree $\frac{-1}{\rho-1},-1$ and
$\rho$, respectively, then $u$ is a self-similar solution, i.e.,
\[
u(x,t)\equiv u_{\lambda}(x,t):=\lambda^{\frac{1}{\rho-1}}u(\lambda
x,\lambda^{2}t),\;\text{ for all }\lambda>0.
\]

\end{theorem}

\subsection{Symmetries and positivity}

In this subsection, we are concerned with symmetry and positivity of
solutions. It is easy to see that the fundamental solution (\ref{Gxyt}) is
positive and invariant by the set $\mathcal{O}_{x_{n}}$ of all rotations
around the axis $\overrightarrow{Ox_{n}}$. Because of that, it is natural to
wonder whether solutions obtained in Theorem \ref{existence} present
positivity and symmetry properties, under certain conditions on the data and potential.

For that matter, let $\mathcal{A}$ be a subset of $\mathcal{O}_{x_{n}}.$ We
recall that a function $f$ is symmetric under the action of $\mathcal{A}$ when
$f(x)=f(T(x))$ for any $T\in\mathcal{A}.$ If $f(x)=-f(T(x))$ for all
$T\in\mathcal{A},$ then $f$ is said to be antisymmetric under $\mathcal{A}.$

\begin{theorem}
\label{pos} Under the hypotheses of Theorem \ref{existence}. Let
$\mathcal{U}\subset\mathbb{R}_{+}^{n}$ be a positive-measure set and
$\mathcal{A}$ a subset of $\mathcal{O}_{x_{n}}$.

\begin{itemize}
\item[(A)] Let $h(a)\geq0$ (resp. $\leq0$) when $a\geq0$ (resp. $\leq0$)$.$ If
$u_{0}\geq0$ (resp. $\leq0$) a.e. in $\mathbb{R}_{+}^{n}$, $u_{0}>0$ (resp.
$<0$) in $\mathcal{U},$ and $V\geq0$ in $\partial\mathbb{R}_{+}^{n},$ then $u$
is positive (resp. negative) in $\mathbb{R}_{+}^{n}\times(0,\infty).$

\item[(B)] Let $h(a)=-h(-a),$ for all $a\in\mathbb{R},$ and let $V$ be
symmetric under the action of $\mathcal{A}|_{\partial\mathbb{R}_{+}^{n}}.$ For
all $t>0,$ the solution $u(\cdot,t)$ is symmetric (resp. antisymmetric), when
$u_{0}$ is symmetric (resp. antisymmetric) under $\mathcal{A}.$
\end{itemize}
\end{theorem}

\begin{remark}
(Special cases of symmetries) Let $h(a)=-h(-a),$ for all $a\in\mathbb{R}.$

\begin{itemize}
\item[(i)] Consider $\mathcal{A}=\mathcal{O}_{x_{n}}$ and let $V$ be radially
symmetric on $\mathbb{R}^{n-1}$. We obtain from item $(B)$ that if $u_{0}$ is
invariant under rotations around the axis $\overrightarrow{Ox_{n}}$ then
$u(\cdot,t)$ does so, for all $t>0$.

\item[(ii)] Let $\mathcal{A}=\{T_{x_{n}}\}$ where $T_{x_{n}}$ is the
reflection with respect to $\overrightarrow{Ox_{n}}$, i.e., $T_{x_{n}%
}((x^{\prime},x_{n}))=(-x^{\prime},x_{n})$ for all $x=(x^{\prime},x_{n})$ and
$x_{n}\geq0$. A function $f$ is said to be $\overrightarrow{Ox_{n}}$-even
(resp. $\overrightarrow{Ox_{n}}$-odd) when $f$ is symmetric (resp.
antisymmetric) under $\{T_{x_{n}}\}.$ If $V(x)$ is an even function then the
solution $u(\cdot,t)$ is $\overrightarrow{Ox_{n}}$-even (resp.
$\overrightarrow{Ox_{n}}$-odd), for all $t>0$, provided that $u_{0}$ is
$\overrightarrow{Ox_{n}}$-even (resp. $\overrightarrow{Ox_{n}}$-odd).
\end{itemize}
\end{remark}

\begin{remark}
Combining Theorems \ref{self-similar} and \ref{pos}, we can obtain solutions
that are both self-similar and invariant by rotations around $\overrightarrow
{Ox_{n}}.$ For instance, in the case $h(a)=\pm\left\vert a\right\vert
^{\rho-1}a$, just take
\[
V(x^{\prime})=\kappa|x^{\prime}|^{-1}\text{ and }u_{0}(x)=\theta\left(
\frac{x_{n}}{\left\vert x\right\vert }\right)  |x|^{-\frac{1}{\rho-1}}%
\]
where $\theta(z)\in BC(\mathbb{R})$ and $\kappa$ is a constant.
\end{remark}

\section{ Proofs}

\label{proofs}

This section is devoted to the proofs of the results. We start by estimating
in Lorentz spaces some linear operators appearing in the integral formulation
(\ref{int})

\subsection{Linear Estimates}

Let $f|_{0}=f(x^{\prime},0)$ stand for the restriction of $f$ to
$\partial\mathbb{R}_{+}^{n}=\mathbb{R}^{n-1}$. We also denote by
$\{E(t)\}_{t\geq0}$ the heat semigroup in the half-space, namely
\begin{equation}
E(t)f(x)=\int_{\mathbb{R}_{+}^{n}}G(x,y,t)f(y)dy \label{heat-1}%
\end{equation}
where $G(x,y,t)$ is the fundamental solution given in (\ref{Gxyt}). For
$\delta>0$ and $1\leq$ $q_{1}\leq q_{2}\leq\infty,$ let us recall the
well-known $L^{q}$-estimate for the heat semigroup $\{e^{t\Delta}\}_{t\geq0}$
on $\mathbb{R}^{n}$:%
\begin{equation}
\Vert(-\Delta_{x})^{\frac{\delta}{2}}e^{t\Delta}f\Vert_{L^{q_{2}}%
(\mathbb{R}^{n})}\leq Ct^{-\frac{1}{2}\left(  \frac{n}{q_{1}}-\frac{n}{q_{2}%
}\right)  -\frac{\delta}{2}}\Vert f\Vert_{L^{q_{1}}(\mathbb{R}^{n})},
\label{heat-ineq-1}%
\end{equation}
where $C>0$ is a constant independent of $f$ and $t$, and $(-\Delta
_{x})^{\frac{\delta}{2}}$ stands for the Riesz potential. For $0<\delta<n$ and
$1\leq q_{1}<q_{2}<\infty$ such that $\frac{1}{\delta}<q_{1}<\frac{n}{\delta}$
and $\frac{n-1}{q_{2}}=\frac{n}{q_{1}}-\delta,$ we have the Sobolev trace-type
inequality in $L^{p}$ (see \cite[Theorem 2]{Adams1})
\begin{equation}
\left\Vert f|_{0}\right\Vert _{L^{q_{2}}(\partial\mathbb{R}_{+}^{n})}\leq
C\left\Vert (-\Delta_{x})^{\frac{\delta}{2}}f\right\Vert _{L^{q_{1}%
}(\mathbb{R}^{n})}. \label{trace-ineq1}%
\end{equation}

The next lemma provide a boundary estimate for (\ref{heat-1}) in the setting
of Lorentz spaces.

\begin{lemma}
\label{prop-trace}Let $1<d_{1}<d_{2}<\infty$ and $1\leq r\leq\infty.$ Then
there exists a constant $C>0$ such that
\begin{equation}
\Vert\lbrack E(t)f]|_{0}\Vert_{L^{(d_{2},r)}(\partial\mathbb{R}_{+}^{n})}\leq
Ct^{-\left(  \frac{n}{2d_{1}}-\frac{n-1}{2d_{2}}\right)  }\Vert f\Vert
_{L^{(d_{1},r)}(\mathbb{R}_{+}^{n})}, \label{est-1}%
\end{equation}
for all $f\in L^{(d_{1},r)}(\mathbb{R}_{+}^{n})$ and $t>0$.
\end{lemma}

\noindent\textbf{Proof.} Consider the extension from $\mathbb{R}_{+}^{n}$ to
$\mathbb{R}^{n}$
\[
\tilde{f}(x)=%
\begin{cases}
f(x^{\prime},x_{n}), & \;x_{n}>0\\
f(x^{\prime},-x_{n}), & \;x_{n}\leq0.
\end{cases}
\]
Now notice that%
\[
E(t)f(x)=e^{t\Delta}\tilde{f}(x)=(g(\cdot,t)\ast\tilde{f})(x)
\]
where
\begin{equation}
g(x,t)=(4\pi t)^{-{n}/{2}}e^{-|x|^{2}/{4t}} \label{heat-kernel1}%
\end{equation}
is the heat kernel on the whole space $\mathbb{R}^{n}$. Therefore,
$e^{t\Delta}\tilde{f}(x)$ is an extension from $\mathbb{R}_{+}^{n}$ to
$\mathbb{R}^{n}$ of $E(t)f(x)$ and
\begin{equation}
\lbrack E(t)f]|_{0}=[e^{t\Delta}\tilde{f}]|_{0}. \label{aux-semigroup1}%
\end{equation}

Let $0<\delta<1$, $\frac{1}{\delta}<l<\frac{n}{\delta}$ and $d_{2}>r$ be such
that $\frac{n-1}{d_{2}}=\frac{n}{l}-\delta.$ It follows from
(\ref{aux-semigroup1}) and (\ref{trace-ineq1}) that
\begin{align*}
\Vert\lbrack E(t)f]|_{0}\Vert_{L^{d_{2}}(\partial\mathbb{R}_{+}^{n})}  &  \leq
C\Vert(-\Delta_{x})^{\frac{\delta}{2}}e^{t\Delta}\tilde{f}\Vert_{L^{l}%
(\mathbb{R}^{n})}\\
&  \leq Ct^{-\frac{1}{2}\left(  \frac{n}{d_{1}}-\frac{n}{l}\right)
-\frac{\delta}{2}}\Vert\tilde{f}\Vert_{L^{d_{1}}(\mathbb{R}^{n})}\\
&  \leq Ct^{-\frac{1}{2}\left(  \frac{n}{d_{1}}-\frac{n-1}{d_{2}}\right)
}\Vert f\Vert_{L^{d_{1}}(\mathbb{R}_{+}^{n})}.
\end{align*}
Now a real interpolation argument leads us
\[
\Vert\lbrack E(t)f]|_{0}\Vert_{L^{(d_{2},r)}(\partial\mathbb{R}_{+}^{n})}\leq
Ct^{-\left(  \frac{n}{2d_{1}}-\frac{n-1}{2d_{2}}\right)  }\Vert f\Vert
_{L^{(d_{1},r)}(\mathbb{R}_{+}^{n})}.
\]
\fin

\bigskip

Let us define the integral operators
\[
\mathcal{G}_{1}(\varphi)(x,t)=\int_{\partial\mathbb{R}_{+}^{n}}G(x,y^{\prime
},t)\varphi(y^{\prime})dy^{\prime}\;\text{ and }\;\mathcal{G}_{2}%
(\varphi)(y^{\prime},t)=\int_{\mathbb{R}_{+}^{n}}G(x,y^{\prime},t)\varphi
(x)dx
\]
where $G(x,y,t)$ is defined in (\ref{Gxyt}). Notice that the functions
$G(x,y^{\prime},t)$ and $\mathcal{G}_{1}(\varphi)(x,t)$ are also well-defined
for $x=(x^{\prime},x_{n})\in\mathbb{R}^{n}$. Recall the pointwise estimate for
the heat kernel (\ref{heat-kernel1}) on $\mathbb{R}^{n}$
\begin{equation}
|(-\Delta_{x})^{\frac{\delta}{2}}g(x,t)|\leq\frac{C_{\delta}}{(t+|x|^{2}%
)^{\frac{n}{2}+\frac{\delta}{2}}}\;\;(\delta\geq0), \label{point-est}%
\end{equation}
for all $x\in\mathbb{R}^{n}$ and $t>0$.

\begin{lemma}
\label{nonl-key-lemma} Let $1<d_{1}<d_{2}<\infty$ and $1\leq r\leq\infty$.
Then, there exists $C>0$ such that
\begin{align}
\Vert\mathcal{G}_{1}(\psi)(x^{\prime},0,t)\Vert_{L^{(d_{2},r)}(\partial
\mathbb{R}_{+}^{n},dx^{\prime})}  &  \leq Ct^{-\left(  \frac{n-1}{2d_{1}%
}-\frac{n-1}{2d_{2}}+\frac{1}{2}\right)  }\Vert\psi\Vert_{L^{(d_{1}%
,r)}(\partial\mathbb{R}_{+}^{n},dx^{\prime})}\label{g1-boundary}\\
\Vert\mathcal{G}_{1}(\psi)(x,t)\Vert_{L^{(d_{2},r)}(\mathbb{R}_{+}^{n},dx)}
&  \leq Ct^{-\left(  \frac{n-1}{2d_{1}}-\frac{n}{2d_{2}}+\frac{1}{2}\right)
}\Vert\psi\Vert_{L^{(d_{1},r)}(\partial\mathbb{R}_{+}^{n},dx^{\prime})}
\label{g1-boundary2}%
\end{align}
for all $\psi\in L^{(d_{1},r)}(\partial\mathbb{R}_{+}^{n}).$
\end{lemma}

\noindent\textbf{Proof.} Let $0<\delta<1$, $\frac{1}{\delta}<l<\frac{n}%
{\delta}$ and $\frac{n-1}{d_{2}}=\frac{n}{l}-\delta$. Firstly, notice that the
trace-type inequality (\ref{trace-ineq1}) yields
\begin{equation}
\Vert\mathcal{G}_{1}(\varphi)(x^{\prime},0,t)\Vert_{L^{d_{2}}(\partial
\mathbb{R}_{+}^{n})}\leq C\Vert(-\Delta_{x})^{\frac{\delta}{2}}\mathcal{G}%
_{1}(\varphi)(x^{\prime},x_{n},t)\Vert_{L^{l}(\mathbb{R}^{n})}. \label{d2}%
\end{equation}
Next we employ Minkowski's inequality for integrals and the pointwise estimate
(\ref{point-est}) to obtain
\begin{align}
\Vert(-\Delta_{x})^{\frac{\delta}{2}}\mathcal{G}_{1}(\varphi)(x^{\prime}%
,x_{n},t)\Vert_{L^{l}(\mathbb{R},dx_{n})}  &  =\left(  \int_{-\infty}^{\infty
}|(-\Delta_{x})^{\frac{\delta}{2}}\mathcal{G}_{1}(\varphi)(x^{\prime}%
,x_{n},t)|^{l}dx_{n}\right)  ^{\frac{1}{l}}\nonumber\\
&  \leq C\int_{\partial\mathbb{R}_{+}^{n}}\left(  \int_{-\infty}^{\infty}%
\frac{|\varphi(y^{\prime})|^{l}}{(t+|x^{\prime}-y^{\prime}|^{2}+x_{n}%
^{2})^{^{\frac{nl}{2}+\frac{\delta l}{2}}}}dx_{n}\right)  ^{\frac{1}{l}%
}dy^{\prime}\nonumber\\
&  =2C\int_{\partial\mathbb{R}_{+}^{n}}|\varphi(y^{\prime})|\left(  \int
_{0}^{\infty}\frac{dx_{n}}{(t+|x^{\prime}-y^{\prime}|^{2}+x_{n}^{2}%
)^{^{\frac{nl}{2}+\frac{\delta l}{2}}}}\right)  ^{\frac{1}{l}}dy^{\prime
}\nonumber\\
&  =2C\int_{\partial\mathbb{R}_{+}^{n}}|\varphi(y^{\prime})|(t+|x^{\prime
}-y^{\prime}|^{2})^{^{-\frac{n}{2}-\frac{\delta}{2}+\frac{1}{2l}}}dy^{\prime
}\label{g1est1}\\
&  \leq Ct^{-\left(  \frac{n}{2}+\frac{\delta}{2}-\frac{1}{2l}\right)  \theta
}\int_{\partial\mathbb{R}_{+}^{n}}|x^{\prime}-y^{\prime}|^{-2\left(  \frac
{n}{2}+\frac{\delta}{2}-\frac{1}{2l}\right)  (1-\theta)}|\varphi(y^{\prime
})|dy^{\prime} \label{g1est2}%
\end{align}

\noindent where (\ref{g1est2}) is obtained from (\ref{g1est1}) by using that
$(a+b)^{-k}\leq a^{-k\theta}b^{-k(1-\theta)}$ when $0<\theta<1$ and
$\kappa\geq0$. Let $d_{1}<l$ and $\gamma=(n-1)(\frac{1}{d_{1}}-\frac{1}{l})$.
Let $0<\theta<1$ be such that $(n-1)-\gamma=\left(  n+\delta-\frac{1}%
{l}\right)  (1-\theta)$. It follows that $\frac{1}{l}=\frac{1}{d_{1}}%
-\frac{\gamma}{n-1}>0$, and Sobolev embedding theorem gives us
\begin{equation}
\left\Vert \int_{\partial\mathbb{R}_{+}^{n}}\frac{1}{|x^{\prime}-y^{\prime
}|^{(n-1)-\gamma}}|\varphi(y^{\prime})|dy^{\prime}\right\Vert _{L^{l}%
(\partial\mathbb{R}_{+}^{n})}\leq C\Vert\varphi\Vert_{L^{d_{1}}(\partial
\mathbb{R}_{+}^{n})}. \label{aux-d5}%
\end{equation}
Fubini's theorem, (\ref{g1est2}) and (\ref{aux-d5}) imply that
\begin{align}
\Vert(-\Delta_{x})^{\frac{\delta}{2}}\mathcal{G}_{1}(\varphi)(x^{\prime}%
,x_{n},t)\Vert_{L^{l}(\mathbb{R}^{n})}  &  =\left\Vert \Vert(-\Delta
_{x})^{\frac{\delta}{2}}\mathcal{G}_{1}(\varphi)(x^{\prime},x_{n}%
,t)\Vert_{L^{l}(\mathbb{R},dx_{n})}\right\Vert _{L^{l}(\mathbb{R}%
^{n-1},dx^{\prime})}\nonumber\\
&  \leq Ct^{-\left(  \frac{n}{2}+\frac{\delta}{2}-\frac{1}{2l}\right)  \theta
}\left\Vert \int_{\partial\mathbb{R}_{+}^{n}}\frac{1}{|x^{\prime}-y^{\prime
}|^{(n-1)-\gamma}}|\varphi(y^{\prime})|dy^{\prime}\right\Vert _{L^{l}%
(\partial\mathbb{R}_{+}^{n},dx^{\prime})}\nonumber\\
&  \leq Ct^{-\left(  \frac{n}{2}+\frac{\delta}{2}-\frac{1}{2l}\right)  \theta
}\Vert\varphi\Vert_{L^{d_{1}}(\partial\mathbb{R}_{+}^{n})}. \label{d5}%
\end{align}
It follows from (\ref{d5}) and (\ref{d2}) that
\begin{equation}
\Vert\mathcal{G}_{1}(\varphi)(x^{\prime},0,t)\Vert_{L^{d_{2}}(\partial
\mathbb{R}_{+}^{n})}\leq Ct^{-\left(  \frac{n}{2}+\frac{\delta}{2}-\frac
{1}{2l}\right)  \theta}\Vert\varphi\Vert_{L^{d_{1}}(\partial\mathbb{R}_{+}%
^{n})}=Ct^{\frac{n-1}{2d_{2}}-\frac{n-1}{2d_{1}}-\frac{1}{2}}\Vert\varphi
\Vert_{L^{d_{1}}(\partial\mathbb{R}_{+}^{n})}, \label{aux-ineq1}%
\end{equation}
because of the equality
\begin{align*}
-\left(  \frac{n}{2}+\frac{\delta}{2}-\frac{1}{2l}\right)  \theta &
=\frac{n-1}{2}-\frac{\gamma}{2}-\frac{1}{2}\left(  n+\delta-\frac{1}{l}\right)
\\
&  =\frac{n-1}{2d_{2}}-\frac{n-1}{2d_{1}}-\frac{1}{2}.
\end{align*}
Now the estimate (\ref{g1-boundary}) follows from (\ref{aux-ineq1}) and real
interpolation. The proof of (\ref{g1-boundary2}) is similar and is left to the reader.\fin

\bigskip

In the next lemma we obtain refined boundary estimates on the Lorentz space
$L^{(d,1)}$ that is the pre-dual one of $L^{(d^{\prime},\infty)}.$ These can
be seen as extensions of Yamazaki's estimates (see \cite{Yamazaki}) to the
operators $\mathcal{G}_{1}$ and $\mathcal{G}_{2}$.

\begin{lemma}
\label{prop-key1} Let $1<d_{1}<d_{2}<\infty$, there exists a constant $C>0$
such that
\begin{align}
\int_{0}^{\infty}t^{\left(  \frac{n-1}{2d_{1}}-\frac{n-1}{2d_{2}}\right)
-\frac{1}{2}}\Vert\mathcal{G}_{1}(\psi)(\cdot,0,t)\Vert_{L^{(d_{2}%
,1)}(\partial\mathbb{R}_{+}^{n})}dt  &  \leq C\Vert\psi\Vert_{L^{(d_{1}%
,1)}(\partial\mathbb{R}_{+}^{n})}\label{prop-key2-est1}\\
\int_{0}^{\infty}t^{\left(  \frac{n-1}{2d_{1}}-\frac{n}{2d_{2}}\right)
-\frac{1}{2}}\Vert\mathcal{G}_{1}(\psi)(\cdot,t)\Vert_{L^{(d_{2}%
,1)}(\mathbb{R}_{+}^{n})}dt  &  \leq C\Vert\psi\Vert_{L^{(d_{1},1)}%
(\partial\mathbb{R}_{+}^{n})}\label{prop-key2-est2}\\
\int_{0}^{\infty}t^{\left(  \frac{n}{2d_{1}}-\frac{n-1}{2d_{2}}\right)
-1}\Vert\mathcal{G}_{2}(\varphi)(\cdot,t)\Vert_{L^{(d_{2},1)}(\partial
\mathbb{R}_{+}^{n})}dt  &  \leq C\Vert\varphi\Vert_{L^{(d_{1},1)}%
(\mathbb{R}_{+}^{n})}, \label{prop-key1-est1}%
\end{align}
for all $\psi\in L^{(d_{1},1)}(\partial\mathbb{R}_{+}^{n})$ and $\varphi\in
L^{(d_{1},1)}(\mathbb{R}_{+}^{n}).$
\end{lemma}

\noindent\textbf{Proof. }We start with (\ref{prop-key1-est1}). Let
$1<p_{1}<d_{1}<p_{2}<d_{2}$ be such that $\frac{1}{d_{1}}-\frac{1}{p_{1}%
}<-\frac{2}{n}$ and $\frac{1}{d_{1}}-\frac{1}{p_{2}}<2$. Noting that
$\mathcal{G}_{2}(\varphi)(y^{\prime},t)=[E(t)\varphi](y^{\prime},0)$, Lemma
\ref{prop-trace} yields
\begin{equation}
\Vert\mathcal{G}_{2}(\varphi)(\cdot,t)\Vert_{L^{(d_{2},1)}(\partial
\mathbb{R}_{+}^{n})}\leq Ct^{-\left(  \frac{n}{2p_{k}}-\frac{n-1}{2d_{2}%
}\right)  }\Vert\varphi\Vert_{L^{(p_{k},1)}(\mathbb{R}_{+}^{n})},\text{ for
}k=1,2. \label{g2-boundary-aux}%
\end{equation}
For $\varphi\in L^{(p_{1},\infty)}(\mathbb{R}_{+}^{n})\cap L^{(p_{2},\infty
)}(\mathbb{R}_{+}^{n})$, we define the following sub-linear operator
\[
\mathcal{F}(\varphi)(t)=t^{\frac{n}{2d_{1}}-\frac{n-1}{2d_{2}}-1}%
\Vert\mathcal{G}_{2}(\varphi)(\cdot,t)\Vert_{L^{(d_{2},1)}(\partial
\mathbb{R}_{+}^{n})}.
\]
Since $1<p_{k}<d_{2}$, it follows from (\ref{g2-boundary-aux}) that
\[
\mathcal{F}(\varphi)(t)\leq Ct^{\left(  \frac{n}{2d_{1}}-\frac{n}{2p_{k}%
}\right)  -1}\Vert\varphi\Vert_{L^{(p_{k},1)}(\mathbb{R}_{+}^{n})}.
\]
Let $\frac{1}{s_{k}}=1-\left(  \frac{n}{2d_{1}}-\frac{n}{2p_{k}}\right)  $ and
take $0<\theta<1$ such that $\frac{1}{d_{1}}=\frac{1-\theta}{p_{1}}%
+\frac{\theta}{p_{2}}$. Then $\frac{1-\theta}{s_{1}}+\frac{\theta}{s_{2}}=1$
with $0<s_{1}<1<s_{2}$. Therefore,
\begin{align*}
\Vert\mathcal{F}(\varphi)(t)\Vert_{L^{(s_{k},\infty)}(0,\infty)}  &  \leq
C\left\Vert t^{-1/s_{k}}\right\Vert _{L^{(s_{k},\infty)}(0,\infty)}^{\ast
}\Vert\varphi\Vert_{L^{(p_{k},1)}(\mathbb{R}_{+}^{n})}\\
&  \leq C\Vert\varphi\Vert_{L^{(p_{k},1)}(\mathbb{R}_{+}^{n})},
\end{align*}
and so $\mathcal{F}:L^{(p_{k},1)}(\mathbb{R}_{+}^{n})\rightarrow
L^{(s_{k},\infty)}(0,\infty)$ is a bounded sublinear operator, for $k=1,2$.
Taking
\[
m_{k}=\Vert\mathcal{F}(\varphi)\Vert_{L^{(p_{k},1)}(\mathbb{R}_{+}%
^{n})\rightarrow L^{(s_{k},\infty)}(0,\infty)},
\]
and recalling the interpolation properties
\[
L^{(d_{1},1)}=(L^{(p_{1},1)},L^{(p_{2},1)})_{\theta,1}\text{ and }%
L^{1}=(L^{(s_{1},\infty)},L^{(s_{2},\infty)})_{\theta,1},
\]
we obtain
\[
\Vert\mathcal{F}(\varphi)\Vert_{L^{1}(0,\infty)}\leq Cm_{1}^{1-\theta}%
m_{2}^{\theta}\Vert\varphi\Vert_{L^{(d_{1},1)}(\mathbb{R}_{+}^{n})}\leq
C\Vert\varphi\Vert_{L^{(d_{1},1)}(\mathbb{R}_{+}^{n})},
\]
which is exactly (\ref{prop-key1-est1}).

In order to show (\ref{prop-key2-est1}), now we define
\[
\mathcal{F}(\psi)(t)=t^{\frac{n-1}{2d_{1}}-\frac{n-1}{2d_{2}}-\frac{1}{2}%
}\Vert\mathcal{G}_{1}(\psi)(\cdot,0,t)\Vert_{L^{(d_{2},1)}(\partial
\mathbb{R}_{+}^{n})}%
\]
and obtain by means of (\ref{g1-boundary}) that
\[
\mathcal{F}(\psi)(t)\leq Ct^{\frac{n-1}{2d_{1}}-\frac{n-1}{2p_{k}}-1}\Vert
\psi\Vert_{L^{(p_{k},1)}(\partial\mathbb{R}_{+}^{n})}.
\]
Let $\frac{1}{s_{k}}=1-\left(  \frac{n-1}{2d_{1}}-\frac{n-1}{2p_{k}}\right)  $
and $0<\theta<1$ be such that $\frac{1}{d_{1}}=\frac{1-\theta}{p_{1}}%
+\frac{\theta}{p_{2}}$. Then $\frac{1-\theta}{s_{1}}+\frac{\theta}{s_{2}}=1$,
and one can obtain (\ref{prop-key2-est1}) by proceeding similarly to proof of
(\ref{prop-key1-est1}). The proof of (\ref{prop-key2-est2}) follows
analogously by considering
\[
\mathcal{F}(\psi)(t)=t^{\frac{n-1}{2d_{1}}-\frac{n}{2d_{2}}-\frac{1}{2}}%
\Vert\mathcal{G}_{1}(\psi)(\cdot,\cdot,t)\Vert_{L^{(d_{2},1)}(\mathbb{R}%
_{+}^{n})}%
\]
and using (\ref{g1-boundary2}) instead of (\ref{g1-boundary}). The details are
left to the reader. \fin

\bigskip

\subsection{Nonlinear estimates}

This section is devoted to estimate the operators
\begin{align}
\mathcal{N}(u)(x,t)  &  =\int_{0}^{t}\int_{\partial\mathbb{R}_{+}^{n}%
}G(x,y^{\prime},t-s)h(u(y^{\prime},s))dy^{\prime}ds\label{opN}\\
\mathcal{T}(u)(x,t)  &  =\int_{0}^{t}\int_{\partial\mathbb{R}_{+}^{n}%
}G(x,y^{\prime},t-s){V}(y^{\prime})u(y^{\prime},s)dy^{\prime}ds. \label{opT}%
\end{align}
For that matter, we define (for each fixed $t>0$)
\[
k_{t}(x,y^{\prime},s)=%
\begin{cases}
G(x,y^{\prime},s), & \text{if \ }0<s<t\\
& \\
0, & \text{otherwise}%
\end{cases}
\]
and consider $\mathcal{H}$ the boundary parabolic integral operator
\[
\mathcal{H}(f)(x,t)=\int_{0}^{\infty}\int_{\partial\mathbb{R}_{+}^{n}}%
k_{t}(x,y^{\prime},t-s)f(y^{\prime},s)dy^{\prime}ds.
\]
For a suitable function $\varphi$ defined in either $\Omega=\mathbb{R}_{+}%
^{n}$ or $\Omega=\partial\mathbb{R}_{+}^{n},$ let us denote
\[
\langle\mathcal{H}(f),\varphi\rangle_{\Omega}=\int_{\Omega}\mathcal{H}%
(f)(x,t)\varphi(x)dx.
\]
Using Tonelli's theorem, we have that
\begin{align}
\left\vert \langle\mathcal{H}(f),\varphi\rangle_{\mathbb{R}_{+}^{n}%
}\right\vert  &  \leq\int_{\mathbb{R}_{+}^{n}}\int_{0}^{\infty}\left(
\int_{\partial\mathbb{R}_{+}^{n}}k_{t}(x,y^{\prime},t-s)\left\vert
f(y^{\prime},s)\right\vert dy^{\prime}\right)  ds\left\vert \varphi
(x)\right\vert dx\nonumber\\
&  =\int_{0}^{\infty}\int_{\partial\mathbb{R}_{+}^{n}}\left\vert f(y^{\prime
},s)\right\vert \left(  \int_{\mathbb{R}_{+}^{n}}k_{t}(x,y^{\prime
},t-s)\left\vert \varphi(x)\right\vert dx\right)  dy^{\prime}ds\nonumber\\
&  =\int_{0}^{\infty}\langle\left\vert f(\cdot,s)\right\vert ,\mathcal{G}%
_{2}(\left\vert \varphi\right\vert )(\cdot,t-s)\rangle_{\partial\mathbb{R}%
_{+}^{n}}ds \label{aux-oper-1}%
\end{align}
and
\begin{align}
\left\vert \langle\mathcal{H}f,\varphi\rangle_{\partial\mathbb{R}_{+}^{n}%
}\right\vert  &  \leq\int_{0}^{\infty}\int_{\partial\mathbb{R}_{+}^{n}%
}\left\vert f(y^{\prime},s)\right\vert \mathcal{G}_{1}(\left\vert
\varphi\right\vert )(y^{\prime},0,t-s)dy^{\prime}ds\nonumber\\
&  =\int_{0}^{\infty}\langle\left\vert f(\cdot,s)\right\vert ,\mathcal{G}%
_{1}(\left\vert \varphi\right\vert )(\cdot,0,t-s)\rangle_{\partial
\mathbb{R}_{+}^{n}}ds, \label{eq2-par}%
\end{align}
because the kernel of $\mathcal{G}_{1}(\varphi)(\cdot,0,t-s)$ and
$\mathcal{G}_{2}(\varphi)(\cdot,t-s)$ is $k_{t}(x,y^{\prime},t-s)$. \bigskip

\begin{lemma}
\label{key-1}Let $n\geq3,$ $\frac{n-1}{n-2}<\rho<\infty,$ and $q=(n-1)(\rho
-1),$ $p=n(\rho-1)$. There exists a constant $C>0$ such that
\begin{align}
\sup_{t>0}\Vert\mathcal{H}\left(  f\right)  (\cdot,t) \Vert_{L^{(p,\infty
)}(\mathbb{R}_{+}^{n})}  &  \leq C\sup_{t>0}\Vert f(\cdot,t)\Vert
_{L^{(\frac{{q}}{\rho},\infty)}(\partial\mathbb{R}_{+}^{n})}\label{est-key2}\\
\sup_{t>0}\Vert\mathcal{H}\left(  f\right)  (\cdot,t) \Vert_{L^{(q,\infty
)}(\partial\mathbb{R}_{+}^{n})}  &  \leq C\sup_{t>0}\Vert f(\cdot
,t)\Vert_{L^{(\frac{{q}}{\rho},\infty)}(\partial\mathbb{R}_{+}^{n})}
\label{est-key3}%
\end{align}
for all $f\in L^{\infty}((0,\infty);L^{(\frac{q}{\rho},\infty)}(\partial
\mathbb{R}_{+}^{n})).$
\end{lemma}

\noindent\textbf{Proof.} Estimate (\ref{aux-oper-1}) and H\"{o}lder inequality
(\ref{holder}) yields
\begin{align}
|\langle\mathcal{H}f,\varphi\rangle_{\mathbb{R}_{+}^{n}}|  &  \leq C\int
_{0}^{\infty}\Vert f(\cdot,s)\Vert_{L^{(\frac{q}{\rho},\infty)}(\partial
\mathbb{R}_{+}^{n})}\Vert\mathcal{G}_{2}(\left\vert \varphi\right\vert
)(\cdot,t-s)\Vert_{L^{(\frac{q}{q-\rho},1)}(\partial\mathbb{R}_{+}^{n}%
)}ds\nonumber\\
&  \leq C\sup_{t>0}\Vert f(\cdot,t)\Vert_{L^{(\frac{q}{\rho},\infty)}%
(\partial\mathbb{R}_{+}^{n})}\int_{0}^{\infty}\Vert\mathcal{G}_{2}(\left\vert
\varphi\right\vert )(\cdot,t-s)\Vert_{L^{(\frac{q}{q-\rho},1)}(\partial
\mathbb{R}_{+}^{n})}ds. \label{aux-est-H1}%
\end{align}
Next, notice that $(\frac{q}{\rho})^{\prime}=\frac{q}{q-\rho}>p^{\prime}$ and
\[
\frac{1}{2}\left(  \frac{n}{p^{\prime}}-\frac{n-1}{(\frac{q}{\rho})^{\prime}%
}\right)  -1=\frac{1}{2}\left(  1-\frac{1}{\rho-1}+\frac{\rho}{\rho-1}\right)
-1=0.
\]
In view of (\ref{aux-est-H1}), we can use duality and estimate
(\ref{prop-key1-est1}) with $d_{1}=p^{\prime}$ and $d_{2}=\frac{q}{q-\rho}$ to
obtain
\begin{align}
I_{1}(t)  &  =\Vert\mathcal{H}(f)(\cdot,t)\Vert_{L^{(p,\infty)}(\mathbb{R}%
_{+}^{n})}=\sup_{\Vert\varphi\Vert_{L^{(p^{\prime},1)}(\mathbb{R}_{+}^{n})}%
=1}\left\vert \langle\mathcal{H}(f),\varphi\rangle_{\mathbb{R}_{+}^{n}%
}\right\vert \nonumber\\
&  \leq C\sup_{t>0}\Vert f(\cdot,t)\Vert_{L^{(\frac{q}{\rho},\infty)}%
(\partial\mathbb{R}_{+}^{n})}\sup_{\Vert\varphi\Vert_{L^{(p^{\prime}%
,1)}(\mathbb{R}_{+}^{n})}=1}\left(  \int_{0}^{\infty}\Vert\mathcal{G}%
_{2}(\left\vert \varphi\right\vert )(\cdot,t-s)\Vert_{L^{(\frac{q}{q-\rho}%
,1)}(\partial\mathbb{R}_{+}^{n})}ds\right) \nonumber\\
&  \leq C\,\sup_{t>0}\Vert f(\cdot,t)\Vert_{L^{(\frac{q}{\rho},\infty
)}(\partial\mathbb{R}_{+}^{n})}\sup_{\Vert\varphi\Vert_{L^{(p^{\prime}%
,1)}(\mathbb{R}_{+}^{n})}=1}\Vert\varphi\Vert_{L^{(p^{\prime},1)}%
(\mathbb{R}_{+}^{n})}\nonumber\\
&  \leq C\,\sup_{t>0}\Vert f(\cdot,t)\Vert_{L^{(\frac{q}{\rho},\infty
)}(\partial\mathbb{R}_{+}^{n})}, \label{aux-proof-1}%
\end{align}
for a.e. $t>0.$ The estimate (\ref{est-key2}) follows by taking the essential
supremum over $(0,\infty)$ in both sides of (\ref{aux-proof-1}).

Now we deal with (\ref{est-key3}) which is the boundary part of the norm
$\left\Vert \cdot\right\Vert _{\mathcal{X}_{p,q}}$. We have that $\left(
\frac{q}{\rho}\right)  ^{\prime}>q^{\prime}$ and
\[
\frac{1}{2}\left(  \frac{n-1}{q^{\prime}}-\frac{1}{(\frac{q}{\rho})^{\prime}%
}\right)  -\frac{1}{2}=\frac{n-1}{2}\left(  \frac{\rho}{q}-\frac{1}{q}\right)
-\frac{1}{2}=0.
\]
Proceeding similarly to proof of (\ref{aux-proof-1}), but using
(\ref{prop-key2-est1}) instead of (\ref{prop-key1-est1}), we obtain
\begin{align}
I_{2}(t)  &  =\Vert\mathcal{H}(f)(\cdot,t)\Vert_{L^{(q,\infty)}(\partial
\mathbb{R}_{+}^{n})}\nonumber\\
&  \leq\sup_{\Vert\varphi\Vert_{L^{(q^{\prime},1)}(\partial\mathbb{R}_{+}%
^{n})}=1}\int_{0}^{\infty}\Vert f(\cdot,s)\Vert_{L^{(\frac{q}{\rho},\infty
)}(\partial\mathbb{R}_{+}^{n})}\Vert\mathcal{G}_{1}(\left\vert \varphi
\right\vert )(\cdot,0,t-s)\Vert_{L^{(\frac{q}{q-\rho},1)}(\partial
\mathbb{R}_{+}^{n})}ds\nonumber\\
&  \leq C\sup_{t>0}\Vert f(\cdot,t)\Vert_{L^{(\frac{q}{\rho},\infty)}%
(\partial\mathbb{R}_{+}^{n})}\int_{0}^{\infty}\Vert\mathcal{G}_{1}(\left\vert
\varphi\right\vert )(\cdot,0,t-s)\Vert_{L^{(\frac{q}{q-\rho},1)}%
(\partial\mathbb{R}_{+}^{n})}ds\nonumber\\
&  \leq C\sup_{t>0}\Vert f(\cdot,t)\Vert_{L^{(\frac{q}{\rho},\infty)}%
(\partial\mathbb{R}_{+}^{n})}\text{ }\sup_{\Vert\varphi\Vert_{L^{(q^{\prime
},1)}(\partial\mathbb{R}_{+}^{n})}=1}\Vert\varphi\Vert_{L^{(q^{\prime}%
,1)}(\partial\mathbb{R}_{+}^{n})}\nonumber\\
&  =C\sup_{t>0}\Vert f(\cdot,t)\Vert_{L^{(\frac{q}{\rho},\infty)}%
(\partial\mathbb{R}_{+}^{n})}, \label{aux-proof-2}%
\end{align}
for a.e. $t\in(0,\infty),$ which is equivalent to (\ref{est-key3}).\fin

\subsection{Proof of Theorem \ref{existence}}

\ \vspace{-0.4cm}

{\textbf{Part (A)}}: Let us write (\ref{int}) as
\begin{equation}
u=E(t)u_{0}+\mathcal{N}(u)+\mathcal{T}(u)\nonumber
\end{equation}
where the operators $\mathcal{N}$ and $\mathcal{T}$ are defined in (\ref{opN})
and (\ref{opT}), respectively.

Recall the heat estimate (see e.g. \cite[Lemma 3.4]{Ukai})
\begin{equation}
\Vert E(t)u_{0}\Vert_{L^{d_{2}}(\mathbb{R}_{+}^{n})}\leq Ct^{-\frac{n}%
{2}(\frac{1}{d_{2}}-\frac{1}{d_{1}})}\Vert u_{0}\Vert_{L^{d_{1}}%
(\mathbb{R}_{+}^{n})}, \label{aux-est}%
\end{equation}
for $1\leq d_{1}\leq d_{2}\leq\infty.$ By using interpolation, (\ref{aux-est})
leads us to
\begin{equation}
\Vert E(t)u_{0}\Vert_{L^{(d_{2},\infty)}(\mathbb{R}_{+}^{n})}\leq
Ct^{-\frac{n}{2}(\frac{1}{d_{2}}-\frac{1}{d_{1}})}\Vert u_{0}\Vert
_{L^{(d_{1},\infty)}(\mathbb{R}_{+}^{n})}, \label{aux-est2}%
\end{equation}
for $1<d_{1}\leq d_{2}<\infty.$

We consider the Banach space $E=BC((0,\infty);\mathcal{X}_{p,q})$ endowed with
the norm (\ref{norm1}). Estimate (\ref{aux-est2}) and Lemma \ref{prop-trace}
yield%
\begin{align}
\Vert E(t)u_{0}\Vert_{E}  &  =\sup_{t>0}\left\Vert E(t)u_{0}\right\Vert
_{L^{(p,\infty)}(\mathbb{R}_{+}^{n})}+\sup_{t>0}\left\Vert E(t)u_{0}%
\right\Vert _{L^{(q,\infty)}(\partial\mathbb{R}_{+}^{n})}\nonumber\\
&  \leq C\left(  \Vert u_{0}\Vert_{L^{(p,\infty)}(\mathbb{R}_{+}^{n})}+\Vert
u_{0}\Vert_{L^{(p,\infty)}(\mathbb{R}_{+}^{n})}\right) \nonumber\\
&  =\delta_{2}\Vert u_{0}\Vert_{L^{(p,\infty)}(\mathbb{R}_{+}^{n})}%
\leq\varepsilon, \label{linear-aux-1}%
\end{align}
provided that $\Vert u_{0}\Vert_{L^{(p,\infty)}(\mathbb{R}_{+}^{n})}\leq
\frac{\varepsilon}{\delta_{2}}$. In what follows, we estimate the operators
$\mathcal{T}$ and $\mathcal{N}$ in order to employ a contraction argument in
$E$. Since $\frac{\rho}{q}=\frac{1}{q}+\frac{\rho-1}{q}$, property
(\ref{nonli-1}) and H\"{o}lder's inequality (\ref{holder}) yield
\begin{align}
\left\Vert h(u)-h(v)\right\Vert _{L^{(q/\rho,\infty)}(\partial\mathbb{R}%
_{+}^{n})}  &  \leq\eta\Vert\left\vert u-v\right\vert (\left\vert u\right\vert
^{\rho-1}+\left\vert v\right\vert ^{\rho-1})\Vert_{L^{(q/\rho,\infty
)}(\partial\mathbb{R}_{+}^{n})}\nonumber\\
&  \leq C\Vert u-v\Vert_{L^{(q,\infty)}(\partial\mathbb{R}_{+}^{n})}(\Vert
u\Vert_{L^{(q,\infty)}(\partial\mathbb{R}_{+}^{n})}^{\rho-1}+\Vert
v\Vert_{L^{(q,\infty)}(\partial\mathbb{R}_{+}^{n})}^{\rho-1}).
\label{ineq-1-aux}%
\end{align}
Using Lemma \ref{key-1} and (\ref{ineq-1-aux}), we obtain
\begin{align*}
\sup_{t>0}\Vert\mathcal{N}(u)-\mathcal{N}(v)\Vert_{\mathcal{X}_{p,q}}  &
=\sup_{t>0}\Vert\mathcal{H}(h(u)-h(v))\Vert_{\mathcal{X}_{p,q}}\\
&  \leq C\sup_{t>0}\Vert h(u)-h(v)\Vert_{L^{(q/\rho,\infty)}(\partial
\mathbb{R}_{+}^{n})}\\
&  \leq K\Vert u-v\Vert_{E}(\Vert u\Vert_{E}^{\rho-1}+\Vert v\Vert_{E}%
^{\rho-1}).
\end{align*}
Also, noting that $\frac{\rho}{q}=\frac{1}{n-1}+\frac{1}{q},$ we have that
\begin{align*}
\Vert\mathcal{T}(u)-\mathcal{T}(v)\Vert_{E}  &  =\sup_{t>0}\Vert
\mathcal{H}(V(u-v))\Vert_{\mathcal{X}_{p,q}}\\
&  \leq C\sup_{t>0}\Vert V(u-v)\Vert_{L^{(q/\rho,\infty)}(\partial
\mathbb{R}_{+}^{n})}\\
&  \leq\delta_{1}\Vert V\Vert_{L^{(n-1,\infty)}(\partial\mathbb{R}_{+}^{n}%
)}\sup_{t>0}\Vert u(\cdot,t)-v(\cdot,t)\Vert_{L^{(q,\infty)}(\partial
\mathbb{R}_{+}^{n})}\\
&  \leq\gamma\Vert u-v\Vert_{E}\text{ with }0<\gamma<1,
\end{align*}
provided that $\gamma=\delta_{1}\Vert V\Vert_{L^{(n-1,\infty)}(\partial
\mathbb{R}_{+}^{n})}$. Now consider
\begin{equation}
\Phi(u)=E(t)u_{0}+\mathcal{N}(u)+\mathcal{T}(u) \label{map-1}%
\end{equation}
and the closed ball $B_{\varepsilon}=\{u\in\mathcal{X}_{p,q}\,;\,\Vert
u\Vert_{\mathcal{X}_{p,q}}\leq\frac{2\varepsilon}{1-\gamma}\}$ where
$\varepsilon>0$ is chosen in such a way that
\begin{equation}
\left(  \frac{2^{\rho}\varepsilon^{\rho-1}K}{(1-\gamma)^{\rho-1}}%
+\gamma\right)  <1. \label{cond-aux-1}%
\end{equation}
For all $u,v\in B_{\varepsilon}$, we obtain that
\begin{align}
\Vert\Phi(u)-\Phi(v)\Vert_{E}  &  \leq\Vert\mathcal{N}(u)-\mathcal{N}%
(v)\Vert_{E}+\Vert\mathcal{T}(u)-\mathcal{T}(v)\Vert_{E}\nonumber\\
&  \leq\Vert u-v\Vert_{E}(K\Vert u\Vert_{E}^{\rho-1}+K\Vert v\Vert_{E}%
^{\rho-1}+\gamma)\label{phi2}\\
&  \leq\left(  K\frac{2^{\rho}\varepsilon^{\rho-1}}{(1-\gamma)^{\rho-1}%
}+\gamma\right)  \Vert u-v\Vert_{\mathcal{X}_{p,q}}.\nonumber
\end{align}
Noting that $\Phi(0)=E(t)u_{0},$ the estimates (\ref{linear-aux-1}) and
(\ref{phi2}) yield
\begin{align}
\Vert\Phi(u)\Vert_{E}  &  \leq\Vert E(t)u_{0}\Vert_{E}+\Vert\Phi
(u)-\Phi(0)\Vert_{E}\nonumber\\
&  \leq\varepsilon+(K\Vert u\Vert_{E}^{\rho}+\gamma\left\Vert u\right\Vert
_{E})\nonumber\\
&  \leq\varepsilon+\left(  K\frac{2^{\rho}\varepsilon^{\rho}}{(1-\gamma
)^{\rho}}+\gamma\frac{2\varepsilon}{1-\gamma}\right)  \leq\frac{2\varepsilon
}{1-\gamma},\nonumber
\end{align}
for all $u\in B_{\varepsilon}$, because of (\ref{cond-aux-1}). Then the map
$\Phi:B_{\varepsilon}\rightarrow B_{\varepsilon}$ is a contraction and Banach
fixed point theorem assures that there is a unique solution $u\in
B_{\varepsilon}$ for (\ref{int}).

The weak convergence to the initial data as $t\rightarrow0^{+}$ follows from
standard arguments and is left to the reader (see e.g. \cite[Lemma
3.8]{Fer-Vil2}, \cite[Lemmas 3.3 and 4.8]{Kozono-Yamazaki}).

\bigskip

{\textbf{Part (B)}}: Let $u,\tilde{u}\in B_{\varepsilon}$ be two solutions
obtained in item \textbf{(A)} corresponding to pairs $(V,u_{0})$ and
$(\tilde{V},\tilde{u}_{0})$, respectively. We have that
\begin{align*}
\Vert u-\tilde{u}\Vert_{E}  &  \leq\Vert E(t)(u_{0}-\tilde{u}_{0})\Vert
_{E}+\Vert\mathcal{N}(u)-\mathcal{N}(\tilde{u})\Vert_{E}+\Vert\mathcal{T}%
(u)-\mathcal{T}(\tilde{u})\Vert_{E}\\
&  \leq\delta_{2}\Vert u_{0}-\tilde{u}_{0}\Vert_{L^{(p,\infty)}(\mathbb{R}%
_{+}^{n})}+K\Vert u-\tilde{u}\Vert_{E}(\Vert u\Vert_{E}^{\rho-1}+\Vert
\tilde{u}\Vert_{E}^{\rho-1})\\
&  +\Vert\mathcal{H}[(V-\tilde{V})\tilde{u}+V(u-\tilde{u})]\Vert_{E}\\
&  \leq\delta_{2}\Vert u_{0}-\tilde{u}_{0}\Vert_{L^{(p,\infty)}(\mathbb{R}%
_{+}^{n})}+\Vert u-\tilde{u}\Vert_{E}\left(  \frac{2^{\rho}\varepsilon
^{\rho-1}K}{(1-\gamma)^{\rho-1}}\right) \\
&  +\delta_{1}\left(  \Vert V-\tilde{V}\Vert_{L^{(n-1,\infty)}}\left\Vert
\tilde{u}\right\Vert _{E}+\Vert V\Vert_{L^{(n-1,\infty)}}\Vert u-\tilde
{u}\Vert_{E}\right) \\
&  \leq\delta_{2}\Vert u_{0}-\tilde{u}_{0}\Vert_{L^{(p,\infty)}(\mathbb{R}%
_{+}^{n})}+\left(  \frac{2^{\rho}\varepsilon^{\rho-1}K}{(1-\gamma)^{\rho-1}%
}+\gamma\right)  \Vert u-\tilde{u}\Vert_{E}+\frac{2\delta_{1}\varepsilon
}{1-\gamma}\Vert V-\tilde{V}\Vert_{L^{(n-1,\infty)}},
\end{align*}
which gives the desired continuity because of (\ref{cond-aux-1}). \fin

\subsection{Proof of Theorem \ref{self-similar}}

From the fixed point argument in the proof of Theorem \ref{existence}, the
solution $u$ is the limit in the space $E$ of the Picard sequence
\begin{equation}
u_{1}=E(t)u_{0},\;u_{k+1}=u_{1}+\mathcal{N}(u_{k})+\mathcal{T}(u_{k}%
),\;k\in\mathbb{N}, \label{Picard}%
\end{equation}
where $\mathcal{N}$ and $\mathcal{T}$ are defined in (\ref{opN}) and
(\ref{opT}), respectively. Since $u_{0}\in L^{(n(\rho-1),\infty)}%
(\mathbb{R}_{+}^{n})$ and $V\in L^{(n-1,\infty)}(\mathbb{R}^{n-1})$, we can
take $u_{0}$ and $V$ as homogeneous functions of degree $-\frac{1}{\rho-1}$
and $-1$, respectively. Using the kernel property
\begin{equation}
G(x,y,t)=\lambda^{n}G(\lambda x,\lambda y,\lambda^{2}t) \label{scal-kernel-1}%
\end{equation}
and homogeneity of $u_{0}$, we have that
\begin{align*}
u_{1}(\lambda x,\lambda^{2}t)  &  =\int_{\mathbb{R}_{+}^{n}}G(\lambda
x,y,\lambda^{2}t)u_{0}(y)dy\\
&  =\int_{\mathbb{R}_{+}^{n}}\lambda^{n}G(\lambda x,\lambda y,\lambda
^{2}t)u_{0}(\lambda y)dy\\
&  =\lambda^{-\frac{1}{\rho-1}}\int_{\mathbb{R}_{+}^{n}}G(x,y,t)u_{0}%
(y)dy=\lambda^{-\frac{1}{\rho-1}}u_{1}(x,t),
\end{align*}
and then $u_{1}$ is invariant by (\ref{sc1}). Recalling that $f(\lambda
a)=\lambda^{\rho}f(a)$ and assuming that
\[
u_{k}(x,t)=u_{k,\lambda}(x,t):=\lambda^{\frac{1}{\rho-1}}u_{k}(\lambda
x,\lambda^{2}t),\text{ for }k\in\mathbb{N}\text{,}%
\]
we obtain
\begin{align*}
\mathcal{N}(u_{k})(\lambda x,\lambda^{2}t)  &  =\int_{0}^{\lambda^{2}t}%
\int_{\partial\mathbb{R}_{+}^{n}}G(\lambda x,y^{\prime},\lambda^{2}%
t-s)h(u_{k}(y^{\prime},s))dy^{\prime}ds\\
&  =\lambda^{n-1+2}\int_{0}^{t}\int_{\partial\mathbb{R}_{+}^{n}}G(\lambda
x,\lambda y^{\prime},\lambda^{2}(t-s))h(\lambda^{-\frac{1}{\rho-1}}%
\lambda^{\frac{1}{\rho-1}}u_{k}(\lambda y^{\prime},\lambda^{2}s))dy^{\prime
}ds\\
&  =\lambda^{n+1}\int_{0}^{t}\int_{\partial\mathbb{R}_{+}^{n}}\lambda
^{-n}G(x,y^{\prime},t-s)\lambda^{-\frac{\rho}{\rho-1}}h(u_{k}(y^{\prime
},s))dy^{\prime}ds\\
&  =\lambda^{-\frac{1}{\rho-1}}\mathcal{N}(u_{k})(x,t)
\end{align*}
and, similarly, $\mathcal{T}(u_{k})(\lambda x,\lambda^{2}t)=\lambda^{-\frac
{1}{\rho-1}}\mathcal{T}(u_{k})(x,t)$. It follows that%

\[
\lambda^{\frac{1}{\rho-1}}u_{k+1}(\lambda x,\lambda^{2}t)=u_{1}%
(x,t)+\mathcal{N}(u_{k})+\mathcal{T}(u_{k})=u_{k+1}(x,t)
\]
and then, by induction, $u_{k}$ is invariant by (\ref{sc1}) for all
$k\in\mathbb{N}$.

Since the norm $\Vert\cdot\Vert_{E}$ is invariant by (\ref{sc1}) and
$u_{k}\rightarrow u$ in $E$, it is easy to see that $u$ is also invariant by
(\ref{sc1}), that is, it is self-similar.\fin

\subsection{Proof of Theorem \ref{pos}}

\textbf{Part (A):} Let $u_{0}\geq0$ a.e. in $\mathbb{R}_{+}^{n}$ and
$\mathcal{U}\subset\mathbb{R}_{+}^{n}$ be a positive measure set with
$u_{0}>0$ in $\mathcal{U}$. It follows from (\ref{Gxyt}) that
\[
u_{1}(x,t)=\int_{\mathbb{R}_{+}^{n}}G(x,y,t)u_{0}(y)dy>0\text{ \ in }%
\overline{\mathbb{R}_{+}^{n}}\times(0,\infty).
\]
By using that $V$ is nonnegative in $\mathbb{R}^{n-1}$ and $h(a)\geq0$ when
$a\geq0,$ one can see that $\mathcal{N}(u)+\mathcal{T}(u)$ is nonnegative in
$\overline{\mathbb{R}_{+}^{n}}\times(0,\infty)$ provided that $u|_{\partial
\mathbb{R}_{+}^{n}}\geq0$. Then, an induction argument applied to the sequence
(\ref{Picard}) shows that $u_{k}>0$ in $\overline{\mathbb{R}_{+}^{n}}%
\times(0,\infty)$, for all $k\in\mathbb{N}$. Since the convergence in the
space $E$ implies convergence in $L^{(p,\infty)}(\mathbb{R}_{+}^{n})$ and in
$L^{(q,\infty)}(\partial\mathbb{R}_{+}^{n})$ for each $t>0,$ we have that (up
to a subsequence) $u_{k}(\cdot,t)\rightarrow u(\cdot,t)$ a.e. in
$(\mathbb{R}_{+}^{n},dx)$ and a.e. in $(\partial\mathbb{R}_{+}^{n},dx^{\prime
})$ for each $t>0$. It follows that $u$ is a nonnegative function because
pointwise convergence preserves nonnegativity. Since $u_{1}>0,$ then
$u=u_{1}+\mathcal{N}(u)+\mathcal{T}(u)\geq u_{1}+0>0$ in $\overline
{\mathbb{R}_{+}^{n}}\times(0,\infty),$ as desired. The proof of the statement
concerning negativity is left to the reader.

\bigskip

\textbf{Part (B):} We only will prove the antisymmetric part of the statement,
because the symmetric one is analogous. Given a $T\in\mathcal{G},$ we have
\begin{align}
u_{1}(T(x),t)  &  =\int_{\mathbb{R}_{+}^{n}}G(T(x),y,t)u_{0}(y)dy\nonumber\\
&  =\int_{\mathbb{R}_{+}^{n}}\dfrac{1}{(4\pi t)^{\frac{n}{2}}}\left[
e^{-\frac{|T(x)-y|^{2}}{4t}}+e^{-\frac{|T(x)-y^{\ast}|^{2}}{4t}}\right]
u_{0}(y)dy\nonumber\\
&  =\int_{\mathbb{R}_{+}^{n}}\dfrac{1}{(4\pi t)^{\frac{n}{2}}}\left[
e^{-\frac{|T((x-T^{-1}(y))|^{2}}{4t}}+e^{-\frac{|T(x-T^{-1}(y^{\ast}))|^{2}%
}{4t}}\right]  u_{0}(y)dy\nonumber\\
&  =\int_{\mathbb{R}_{+}^{n}}\dfrac{1}{(4\pi t)^{\frac{n}{2}}}\left[
e^{-\frac{|x-T^{-1}(y)|^{2}}{4t}}+e^{-\frac{|x-(T^{-1}(y))^{\ast}|^{2}}{4t}%
}\right]  u_{0}(y)dy\nonumber\\
&  =\int_{\mathbb{R}_{+}^{n}}G(x,T^{-1}(y),t)u(y)dy.\nonumber
\end{align}
Making the change of variable $z=T^{-1}(y)$ and using that $u_{0}$ is
antisymmetric under $\mathcal{G}$, we obtain
\[
u_{1}(T(x),t)=\int_{\mathbb{R}_{+}^{n}}G(x,z,t)u_{0}(T(z))dz=-\int
_{\mathbb{R}_{+}^{n}}G(x,z,t)u_{0}(z)dz=-u_{1}(x,t).
\]
A similar argument shows that
\[
\mathcal{L}(\theta)(x,t)=\int_{0}^{t}\int_{\partial\mathbb{R}_{+}^{n}%
}G(x,y^{\prime},t-s)\theta(y^{\prime},t)dy^{\prime}ds
\]
is antisymmetric when $\theta(\cdot,t)|_{\partial\mathbb{R}_{+}^{n}}$ is also,
for each $t>0$. As $V$ is symmetric and $h(a)=-h(-a)$, it follows that
\[
\theta(x,t)=h(u(\cdot,t))+Vu(\cdot,t)
\]
is antisymmetric whenever $u(\cdot,t)$ does so. Therefore, by means of an
induction argument, one can prove that each element $u_{k}(\cdot,t)$ of the
sequence (\ref{Picard}) is antisymmetric. Recall from Part (A) that (up a
subsequence) $u_{k}(\cdot,t)\rightarrow u(\cdot,t)$ a.e. in $(\mathbb{R}%
_{+}^{n},dx)$ and in $(\partial\mathbb{R}_{+}^{n},dx^{\prime})$, for each
$t>0.$ Since this convergence preserves antisymmetry, it follows that
$u(\cdot,t)$ is also antisymmetric, for each $t>0$. \fin

\end{document}